\documentclass[12pt, a4paper, twoside]{article}
\usepackage{amssymb}
\usepackage[leqno]{amsmath}
\usepackage{latexsym}
\usepackage{enumerate}
\usepackage{amsmath, amssymb, amscd}
\DeclareMathOperator{\rk}{rk}
\DeclareMathOperator{\Hom}{Hom}
\DeclareMathOperator{\Ext}{Ext}

\Roman{equation}

\begin{document}

\title{\textbf{\Large{ACM bundles on K3 surfaces of genus 2}}}

\author{Kenta Watanabe \thanks{Department of Mathematical Sciences, Osaka University, 1-1 Machikaneyama-chou Toyonaka Osaka 560-0043 Japan, {\it E-mail address:goo314kenta@mail.goo.ne.jp}, Telephone numbers: 090-9777-1974} }

\date{}

\maketitle


\begin{abstract}

Let $\pi:X\rightarrow \mathbb{P}^2$ be a K3 surface of genus 2 and $L=\pi^{\ast}\mathcal{O}_{\mathbb{P}^2}(3)$, and assume that $\pi^{\ast}\mathcal{O}_{\mathbb{P}^2}(1)$ is ample as a line bundle on $X$. In this paper, we give a numerical characterization of initialized and ACM line bundles on $X$ with respect to $L$ and construct families of semistable indecomposable ACM bundles of higher rank, by using extensions of ACM line bundles.

\end{abstract}

\section{Introduction}

We work over the complex number field $\mathbb{C}$. Let $X$ be a smooth projective variety and $\mathcal{O}_X(1)$ be a very ample line bundle on $X$. Then, a vector bundle $\mathcal{E}$ on $X$ is called an {\it{Arithmetically Cohen-Macaulay}} ({\it{ACM}} for short) with respect to $\mathcal{O}_X(1)$ if $H^i(X,\mathcal{E}(l))=0$, for $1\leq i\leq\dim(X)-1$ and $l\in\mathbb{Z}$, where $\mathcal{E}(l)=\mathcal{E}\otimes\mathcal{O}_X(l)$.

Previously, many people have studied indecomposable ACM bundles with respect to a given polarization on smooth projective surfaces. In the case where $X$ is a smooth hypersurface in a projective space, if $X$ is a quadric, then ACM bundles on $X$ were completely classified by Kn${\rm{\ddot{o}}}$rrer [Kn]. If $X$ is a cubic surface, Casanellas and Hartshorne [C-H] have constructed an $n^2+1$- dimensional family of rank $n$ indecomposable ACM bundles on $X$ with Chern classes $c_1=nH$ and $c_2=\frac{1}{2}(3n^2-n)$ for $n\geq2$. Moreover, Faenzi [Fa] gave a precise classification of rank 2 ACM bundles on $X$. On the other hand, if $X$ is a quartic surface, we gave a numerical characterization of ACM and initialized line bundles on $X$ [W2]. Here, a line bundle $\mathcal{L}$ is called {\it{initialized}} if $H^0(X,\mathcal{L})\neq0$ and $H^0(X,\mathcal{L}(-1))=0$ are satisfied. Moreover, in the case where $X$ is not a hypersurface, for example, if $X$ is a DelPezzo surface with the very ample anticanonical line bundle $-K_X$, Joan Pons-Llopis and Fabio Tonini [P-T] have classified ACM line bundles on $X$ with respect to $-K_X$, and have constructed families of indecomposable ACM bundles of higher rank, by using extensions of ACM line bundles on $X$.

In this paper, we consider ACM bundles on a K3 surface of genus 2 (i.e., the pair of a K3 surface and an ample line bundle of sectional genus 2). Let $\pi:X\rightarrow\mathbb{P}^2$ be a double cover branched along a smooth plane sextic, and $L:=\pi^{\ast}\mathcal{O}_{\mathbb{P}^2}(3)$. Then, by the characterization of hyperelliptic linear systems on K3 surfaces ([SD], Theorem 5.2), if $\pi^{\ast}\mathcal{O}_{\mathbb{P}^2}(1)$ is ample, $L$ is very ample. Therefore, first of all, we gave a numerical characterization of ACM line bundles with respect to such a line bundle $L$.

\newtheorem{thm}{Theorem}[section]

\begin{thm} Let $X$ and $L$ be as above, and $D$ be a nonzero effective divisor on $X$. Assume that $L$ is ample. Then the following conditions are equivalent.

\smallskip

\smallskip

\noindent {\rm{(i)}} $\mathcal{O}_X(D)$ is an ACM and initialized line bundle with respect to $L$.

\noindent {\rm{(ii)}} For $H\in |L|$, one of the following cases occurs.

\smallskip

{\rm{(a)}} $D^2=-2$ and $H.D\in\{3,6,9\}$.

{\rm{(b)}} $D^2=0$ and $H.D=9$.

{\rm{(c)}} $D^2=2$, $H.D\in \{6,9,12\}$, and if $H.D=12$, then $|H-D|=\emptyset$.

{\rm{(d)}} $D^2=4$ and $H.D\in \{9,12\}$.

{\rm{(e)}} $D^2=8$ and $H.D\in\{12,15\}$.

{\rm{(f)}} $D^2=10$ and $H.D=15$.

{\rm{(g)}} $D^2=14$ and $H.D=18$.

{\rm{(h)}} $D^2\in\{20,26,32\}$,  $D^2=2H.D-22$,

\noindent and $h^0(\mathcal{O}_X(D-H))=h^1(\mathcal{O}_X(2H-D))=0$.

\end{thm}

\noindent In general, a vector bundle $\mathcal{E}$ on a smooth projective surface is called {\it{an Ulrich bundle}} if the following cohomology groups vanish.
$$H^0(X,\mathcal{E}(-1)),\;H^1(X,\mathcal{E}(-1)),\;H^1(X,\mathcal{E}(-2)),\;H^2(X,\mathcal{E}(-2))$$
\noindent Therefore, in Theorem 1.1, an Ulrich line bundle is characterized as an initialized and ACM line bundle with the largest self intersection number among such line bundles.

On the other hand, for a given K3 surface of genus 2 $\pi:X\rightarrow\mathbb{P}^2$, the linear system $|L|$ given by $L$ defined as above is known as a counterexample to a conjecture of Harris and Mumford that the gonality should be constant among the smooth curves in a linear system on K3 surfaces, and it is called Donagi-Morrison's example. In particular, Ciliberto and Pareschi [C-P] proved that if $L$ is ample, then the associated Donagi-Morrison's example $|L|$ is the only counterexample to the conjecture of Harris and Mumford.

In the present work, we proved that if $\pi:X\rightarrow\mathbb{P}^2$ is a K3 surface of genus 2 of the Picard number $\rho(X)\geq2$ whose N$\rm{\acute{e}}$ron-Severi lattice $S_X$ is 2-elementary, that is, there exists a non-negative integer $a$ such that $S_X^{\ast}/S_X\cong(\mathbb{Z}/2\mathbb{Z})^a$, where $S_X^{\ast}:=\Hom(S_X,\mathbb{Z})$, then the Donagi-Morrison's example associated to the double covering $\pi$ can be ample precisely when $\rho(X)=a=9$ ([W1], Theorem 1.2). In general, it is well known that, by the global Torelli Theorem for K3 surfaces, if the N$\rm{\acute{e}}$ron-Severi lattice of a K3 surface is 2-elementary, then there exists a unique involution $\theta$ (it is called the canonical involution) which acts trivially on it and acts like the multiplication by $(-1)$ on the transcendental lattice, and the concrete description of the fixed locus $X^{\theta}$ of $\theta$ is given by Nikulin [Ni]. In particular, if $\rho(X)=a=9$, $X^{\theta}$ is a smooth genus 2 curve which is ample as a divisor on $X$. In [W1], we also proved that if $\rho(X)=a=9$ and the Donagi-Morrison's example associated to the double covering $\pi$ is ample, then it is given by the divisor $3X^{\theta}$. Therefore, in this case, we constructed families of indecomposable semi-stable ACM bundles of higher rank, by using extensions of some ACM line bundles with respect to $L=\mathcal{O}_X(3X^{\theta})$.
\begin{thm} Let $X$ be a K3 surface whose N$\acute{e}$ron-Severi lattice is 2-elementary, and assume that $\rho(X)=a=9$. Let $L=\mathcal{O}_X(3X^{\theta})$. Then, for any integer $n\geq3$, there exists a family of dimension $\geq n$ of non-isomorphic indecomposable semi-stable ACM bundles of rank $n$ with respect to $L$. \end{thm}
Our plan of this paper is as follows. In Section 2, we recall some known results about line bundles and linear systems on K3 surfaces. In Section 3, we give a numerical characterization of ACM line bundles on K3 surfaces with large self intersection numbers. In Section 4, we give a proof of Theorem 1.1. In Section 5, we recall some classical facts about K3 surfaces whose N$\rm{\acute{e}}$ron-Severi lattices are 2-elementary. In Section 6, we recall some known results about indecomposable vector bundles of higher rank. In Section 7, we give a proof of Theorem 1.2.

\smallskip

\smallskip

\noindent {\bf{Notation and Conventions.}} A surface is a smooth projective surface. Let $X$ be a surface. We denote by $K_X$ the canonical line bundle of $X$. For a divisor $D$ on $X$, we will denote by $|D|$ the linear system defined by $D$. If two divisors $D_1$ and $D_2$ on $X$ satisfy the condition that $|D_1|=|D_2|$, then we will write $D_1\sim D_2$. We will denote by $S_X$ the N$\rm{\acute{e}}$ron-Severi lattice of $X$, and denote by $\rho(X)$ its rank. We denote by $\mathcal{O}_X(1)$ a very ample line bundle that provides a closed embedding in a projective space, and denote by $\mathcal{O}_X(l)$ the line bundle $\mathcal{O}_X(1)^{\otimes l}$. For a vector bundle $\mathcal{E}$ on $X$, we will write $\mathcal{E}\otimes\mathcal{O}_X(l)=\mathcal{E}(l)$. We will say that a vector bundle $\mathcal{E}$ on $X$ is initialized with respect to $\mathcal{O}_X(1)$ if it satisfies the condition that
$$ H^0(X,\mathcal{E}(-1))=0 \text{ but } H^0(X,\mathcal{E})\neq0.$$
\noindent We call a regular surface $X$ a K3 surface if the canonical line bundle of $X$ is trivial.

\section{Linear systems and line bundles on K3 surfaces}

In this section, we recall some basic results about ample line bundles and linear systems on K3 surfaces. First of all, we remark some facts about numerical connected divisors on a surface.

\newtheorem{df}{Definition}[section]

\begin{df} A divisor $D$ on a surface is called {\rm{$m$-connected}} if $D_1.D_2\geq m$, for each effective decomposition $D=D_1+D_2$.\end{df}

\noindent If a divisor $D$ on a surface is 1-connected, then $h^0(\mathcal{O}_{D})=1$ (cf. [B-P-W], Corollary 12.3). Hence, we can easily see that, for a 1-connected divisor $D$ on a K3 surface, we get $h^1(\mathcal{O}_X(D))=0$. Next, we recall a result about the classification of base point free divisors on K3 surfaces.

\newtheorem{prop}{Proposition}[section]

\begin{prop}{\rm{([SD], Proposition 2.7)}} Let $L$ be a numerical effective line bundle on a K3 surface $X$. Then $|L|$ is not base point free if and only if there exists an elliptic curve $F$, a smooth rational curve $\Gamma$ and an integer $k\geq2$ such that $F.\Gamma=1$ and $L\cong\mathcal{O}_X(kF+\Gamma)$. \end{prop}

\begin{prop}{\rm{([SD], Proposition 2.6)}} Let $L$ be a line bundle on a K3 surface $X$ such that $|L|\neq\emptyset$. Assume that $|L|$ has no fixed components. Then one of the following cases occurs.

\smallskip

\smallskip

{\rm{(i)}} $L^2>0$ and the general member of $|L|$ is a smooth irreducible curve of genus $\frac{1}{2}L^2+1$.

{\rm{(ii)}} $L^2=0$ and $L\cong\mathcal{O}_X(kF)$, where $k\geq1$ is an integer and $F$ is a smooth curve of genus one. In this case, $h^1(L)=k-1$. \end{prop} 

\noindent It is well known that, for an irreducible curve $C$ on a K3 surface such that $C^2>0$, $|C|$ is base point free ([SD], Theorem 3.1). Hence, by Proposition 2.2, the following proposition follows.

\begin{prop}{\rm{([SD], Corollary 3.2)}} Let $L$ be a line bundle on a K3 surface. Then $|L|$ has no base points outside its fixed components.\end{prop}

\noindent At the end of this section, we recall some classical results about very ample line bundles on K3 surfaces. It is well known that if an ample linear system on a K3 surface is not very ample, then it is hyperelliptic ([SD]). Hence, by the characterization of hyperelliptic linear systems on K3 surfaces, we have the following assertion.

\begin{prop} {\rm{(cf. [M-M], and [SD], Theorem 5.2)}} Let $L$ be a numerical effective line bundle with $L^2\geq4$ on a K3 surface $X$. Then $L$ is very ample if and only if the following conditions are satisfied.

\smallskip

\smallskip

{\rm{(i)}} There is no irreducible curve $E$ such that $E^2=0$ and $E.L=1$ or 2.

{\rm{(ii)}} There is no irreducible curve $E$ such that $E^2=2$ and $L\cong\mathcal{O}_X(2E)$.

{\rm{(iii)}} There is no irreducible curve $E$ such that $E^2=-2$ and $E.L=0$. \end{prop}

\noindent Note that, by Proposition 2.1 and Proposition 2.4, if $L$ is a very ample line bundle, then $|L|$ is base point free. Hence, the general member of it is a smooth irreducible curve. Moreover, by Proposition 2.2, we have the following fact.

$\;$

\noindent {\bf{Corollary 2.1}} {\it{Let}} $X$ {\it{be a K3 surface,}} $L$ {\it{be an ample line bundle, and let}} $D$ {\it{ be a nonzero effective divisor on X with}} $D^2\geq0$. {\it{Then we have the following results.}}

\smallskip

\smallskip 

(i) {\it{If}} $L^2=2$, {\it{then}} $L=\mathcal{O}_X(D)$ {\it{or}} $L.D\geq3$.

(ii) {\it{If}} $L^2=2$ {\it{and}} $L.D=3$, {\it{then}} $D^2=2$ {\it{or}} $|D|$ {\it{is base point free}}.

(iii) {\it{If }} $L$ {\it{is very ample, then}} $L.D\geq 3$.

$\;$

{\it{Proof}}. We take a member $H\in |L|$. (i) First of all, we note that $H.D\geq2$. In fact, if $H.D=1$, we have $(H-2D)^2\geq-2$ and $H.(H-2D)=0$. In this case, by the ampleness of $L$, we have $H\sim 2D$. However, this is a contradiction. If $H.D=2$, we have $(H-D)^2\geq -2$ and $H.(H-D)=0$. Hence, in this case, we have $H\sim D$, by the ampleness of $L$. 

\smallskip

\smallskip

(ii) Assume that $|D|$ is not base point free, and let $\Delta$ be the base divisor of $|D|$ and $D^{'}=D-\Delta$. Since $D^2\geq0$, we have $D^{'}\neq0$. By the ampleness of $L$, we have $H.\Delta\geq1$, and hence, by the assumption and the proof of (i), we have $H.D^{'}=2$ and $H.\Delta=1$. Hence, $\Delta$ is a $(-2)$-curve, and, by the assertion of (i), we have $D^{'}\sim H$. Therefore, we have
$$D^2=(D^{'}+\Delta)^2=2H.\Delta=2.$$

\smallskip

\smallskip

(iii) Note that, since $L$ is very ample, we have $H^2\geq4$. Let $\Delta$ be the base divisor of $|D|$. Since $D^2\geq0$, we have $D-\Delta\neq 0$. If $(D-\Delta)^2=0$, then, by Proposition 2.2, there exists an elliptic curve $F$ with $D-\Delta \sim kF\;(k\geq1)$. In this case, by Proposition 2.4, we have $H.(D-\Delta)\geq 3k\geq3$. If $(D-\Delta)^2\geq 2$, by Hodge index theorem, we have $(H.(D-\Delta))^2\geq (H^2)(D-\Delta)^2\geq8$. Hence, we have the assertion. $\hfill\square$

$\;$

\noindent {\bf{Remark 2.1}}. Let $X$ be a K3 surface, and let $L$ be an ample line bundle with $L^2=2$ on $X$. Then, $L^{\otimes3}$ is very ample.

\section{ACM line bundles on polarized K3 surfaces}

In this section, we remark some known results about ACM line bundles and give a numerical characterization of ACM line bundles with respect to a given very ample line bundle on a K3 surface. First of all, we recall our previous result about ACM line bundles on quartic hypersurfaces in $\mathbb{P}^3$.

\begin{thm} {\rm{([W2], Theorem 1.1)}} Let $X$ be a smooth quartic hypersurface of $\mathbb{P}^3$, and let $D$ be a nonzero effective divisor on $X$. Then the following conditions are equivalent.

\smallskip

\smallskip

\noindent {\rm{(i)}} $\mathcal{O}_X(D)$ is an ACM and initialized line bundle.

\noindent {\rm{(ii)}} For a hyperplane section $H$ of $X$, one of the following cases occurs.

\smallskip

\smallskip

{\rm{(a)}} $D^2=-2$ and $1\leq H.D\leq 3$.

{\rm{(b)}} $D^2=0$ and $3\leq H.D\leq 4$.

{\rm{(c)}} $D^2=2$ and $H.D=5$.

{\rm{(d)}} $D^2=4,\;H.D=6$ and $|D-H|=|2H-D|=\emptyset.$\end{thm}

\noindent In this section, we want to give a numerical characterization of an ACM line bundle $\mathcal{O}_X(D)$ with respect to a very ample line bundle $L$, by using the method as in Theorem 3.1. However, if $D^2$ is sufficiently small compared with $L^2$, it is difficult to do it. Hence, we will only consider the case where $D^2\geq L^2-4$.

\begin{thm} Let $X$ be a K3 surface, and let $L$ be a very ample line bundle. Let $D$ be a nonzero effective divisor on $X$ with $D^2\geq L^2-4$. Then the following conditions are equivalent.

\smallskip

\smallskip

\noindent {\rm{(i)}} $\mathcal{O}_X(D)$ is an ACM and initialized line bundle with respect to $L$.

\noindent {\rm{(ii)}} For $H\in |L|$, one of the following cases occurs.

\smallskip

\smallskip

{\rm{(a)}} $D^2 = H^2-4$ and $H.D\in\{H^2-1,\;H^2\}.$

{\rm{(b)}} $D^2 = H^2-2$ and $H.D=H^2+1.$

{\rm{(c)}} $D^2\geq H^2,\;D^2=2H.D-H^2-4,\;|D-H|=\emptyset$ and $h^1(\mathcal{O}_X(2H-D))=0.$ \end{thm}

\noindent First of all, in order to prove Theorem 3.2, we prepare the following lemmas.

\newtheorem{lem}{Lemma}[section]

\begin{lem} Let $X$ and $L$ be as in Theorem 3.2, and let $D$ be a nonzero effective divisor. Moreover, let $m\in\mathbb{N}$. Then if $L.D\leq mL^2-1$ and, for any $k\in\mathbb{Z}$ with $0\leq k\leq m$,  $h^1(\mathcal{O}_X(D)\otimes L^{-k})=0$, then $\mathcal{O}_X(D)$ is an ACM line bundle. \end{lem}

{\it{Proof}}. Let $m\geq1$ be an integer satisfying the assumption. Let $n\in\mathbb{N}$ and let $H\in |L|$ be a smooth irreducible curve.  First of all, we have
$$h^1(\mathcal{O}_H(nH+D))=h^0(\mathcal{O}_H(-D-(n-1)H))=0.$$
\noindent By the assumption, since we have $h^1(\mathcal{O}_X(D))=0$, by the exact sequence
$$0\rightarrow\mathcal{O}_X(D+(n-1)H)\rightarrow\mathcal{O}_X(D+nH)\rightarrow\mathcal{O}_H(D+nH)\rightarrow0,$$
\noindent we have $h^1(\mathcal{O}_X(D+nH))=0$, by using induction.

On the other hand, since $L.D\leq mL^2-1$, if $n\geq m$, then we have
$$h^1(\mathcal{O}_H((n+1)H-D))=h^0(\mathcal{O}_H(D-nH))=0.$$
\noindent By the exact sequence
$$0\rightarrow\mathcal{O}_X(nH-D)\rightarrow\mathcal{O}_X((n+1)H-D)\rightarrow\mathcal{O}_H((n+1)H-D)\rightarrow0,$$
\noindent we have
$$h^1(\mathcal{O}_X(D-(n+1)H))=h^1(\mathcal{O}_X((n+1)H-D))=0\;(n\geq m).$$
\noindent Hence, we have the assertion, by using induction.$\hfill\square$

\begin{lem} Let $X$ be a K3 surface, and let $D$ be a divisor on $X$ which is not linearly equivalent to 0, and assume that $|D|\neq\emptyset$. Let $\Delta$ be the base divisor of $|D|$. If $h^1(\mathcal{O}_X(D-\Delta))=0$ and $D^2=(D-\Delta)^2$, then $h^1(\mathcal{O}_X(D))=0$.\end{lem}

{\it{Proof}}. Let $D$ be a nonzero effective divisor satisfying the assumption. Then we note that, since $D^2=(D-\Delta)^2$, the movable part of $|D|$ is not empty. Since $h^1(\mathcal{O}_X(D-\Delta))=0$, we have
$$h^0(\mathcal{O}_X(D))=h^0(\mathcal{O}_X(D-\Delta))=\chi(\mathcal{O}_X(D-\Delta)).$$
\noindent On the other hand, since $D^2=(D-\Delta)^2$, we have $\chi(\mathcal{O}_X(D))=\chi(\mathcal{O}_X(D-\Delta))$. Hence, we have $h^1(\mathcal{O}_X(D))=0$. $\hfill\square$

$\;$

{\it{Proof of Theorem 3.2}}. Let $X,\;D$ and $L$ be as in Theorem 3.2. Let $H\in|L|$ be a smooth curve. If $H^2=4$, then the assertion already proved in Theorem 3.1. Hence, we assume that $H^2\geq6$.

\smallskip

\smallskip

\noindent (i)$\Longrightarrow$(ii) First of all, we consider the case where $D^2=H^2-4$ and the case where $D^2=H^2-2$. If $|H-D|=\emptyset$, then, by the assumption, we have
$$\chi(\mathcal{O}_X(H-D))=0,$$
\noindent and hence, we have the assertion. We consider the case where $|H-D|\neq\emptyset$. Note that since $L$ is ample, we have $H.(H-D)>0$. 

Assume that $D^2=H^2-4$. In this case, by Hodge index theorem, we have
$$(H^2-3)^2<(H^2)(D^2)\leq (H.D)^2\leq (H^2-1)^2,$$
\noindent and hence, we have $H.D=H^2-2$ or $H^2-1$. If $H.D=H^2-2$, then $(H-D)^2=0$ and $H.(H-D)=2$. However, by Corollary 2.1 (iii), this is a contradiction. Hence, we have $H.D=H^2-1$. 

Assume that $D^2=H^2-2$. Since $H.(H-D)\geq1$, we have $(H-D)^2\geq0$. Since $L$ is very ample, by Corollary 2.1 (iii), we have $H.(H-D)\geq3$. However, by Hodge index theorem, we have
$$(H^2-2)^2<(H^2)(D^2)\leq (H.D)^2\leq (H^2-3)^2.$$
\noindent This is a contradiction.

Next, we consider the case where $D^2\geq H^2$. In this case, we note that, $|H-D|=\emptyset$. In fact, if $|H-D|\neq\emptyset$, by the ampleness of $L$, we have $H.(H-D)>0$. Hence, by Hodge index theorem, we have
$$(H^2)(D^2)\leq (H.D)^2<(H^2)^2,$$
\noindent and hence, we have the contradiction $D^2<H^2$. Therefore, by the assumption that $\mathcal{O}_X(D)$ is ACM and initialized, the assertion holds.

\smallskip

\smallskip

\noindent (ii)$\Longrightarrow$(i) We consider the case where $D^2=H^2-4$. Since $H.(D-H)=-1$ or 0, we have $|D-H|=\emptyset$. Hence, in this case, $\mathcal{O}_X(D)$ is initialized. 

Assume that $H.D=H^2-1$. By Lemma 3.1, it is sufficient to show that 
$$h^1(\mathcal{O}_X(D))=h^1(\mathcal{O}_X(H-D))=0.$$
\noindent First of all, in order to show that $h^1(\mathcal{O}_X(D))=0$, we show that $|D|$ is base point free. Assume that $|D|$ is not base point free, and $\Delta$ be the base divisor of $|D|$. Let $D^{'}=D-\Delta$. Then we note that, since $H$ is ample, we have $H.D^{'}\leq H^2-2$. Assume that ${D^{'}}^2=0.$ Then, there exists an elliptic curve $F$ and an integer $k\geq1$ such that $D^{'}\sim kF$. Hence, by Corollary 2.1 (iii), we have $k\leq\frac{1}{3}(H^2-2)$. Since, by Proposition 2.2 (ii), $h^1(\mathcal{O}_X(D^{'}))=k-1$, we have
$$\chi(\mathcal{O}_X(D^{'}))\geq\chi(\mathcal{O}_X(D))-\frac{1}{3}(H^2-5).$$
\noindent Hence, we have 
$${D^{'}}^2\geq \frac{1}{3}(H^2-2)>0.$$
\noindent This is a contradiction. Since ${D^{'}}^2>0$, we have $h^1(\mathcal{O}_X(D^{'}))=0$. Hence, by comparing $\chi(\mathcal{O}_X(D^{'}))$ and $\chi(\mathcal{O}_X(D))$, we have ${D^{'}}^2\geq D^2=H^2-4$. By Hodge index theorem, we have
$$(H^2-3)^2<H^2(H^2-4)\leq (H^2)({D^{'}}^2)\leq (H.D^{'})^2\leq (H^2-2)^2.$$
\noindent Since we have $H.D^{'}=H^2-2$, we also have $(H-D^{'})^2\geq0$ and $H.(H-D^{'})=2$. However, this contradicts to the assumption that $H$ is very ample. Therefore, $|D|$ is base point free, and hence, we have $h^1(\mathcal{O}_X(D))=0$. On the other hand, since $(H-D)^2=-2$ and $H.(H-D)=1$, the member of $|H-D|$ is irreducible. Therefore, we have $h^1(\mathcal{O}_X(H-D))=0$, and hence, $\mathcal{O}_X(D)$ is ACM.

Assume that $H.D=H^2$. Since $H.(D-H)=0$ and $D^2=H^2-4$, we have $|H-D|=\emptyset$. Since $(H-D)^2=-4$, we have $h^1(\mathcal{O}_X(H-D))=-\chi(\mathcal{O}_X(H-D))=0$. By Lemma 3.1, it is sufficient to show that $h^1(\mathcal{O}_X(D))=h^1(\mathcal{O}_X(2H-D))=0.$ 

First of all, we consider the case where $|D|$ is base point free. In this case, by the theorem of Bertini, we have $h^1(\mathcal{O}_X(D))=0$. In order to show that $h^1(\mathcal{O}_X(2H-D))=0$, we show that $|2H-D|$ is base point free. Assume that it is not base point free, and let $\Delta$ be the base divisor of it. Then, since $(2H-D)^2=D^2>0$, the divisor $2H-D-\Delta$ is not linearly equivalent to 0. Hence, we take a nonzero divisor $D^{'}\in|2H-D-\Delta|$. If ${D^{'}}^2=0$, then there exists an elliptic curve $F$ and $k\geq1$ such that $D^{'}\sim kF$. Since $H.D^{'}\leq H^2-1$, by Corollary 2.1, we have $k\leq \frac{1}{3}(H^2-1)$. By the same reason as above, we have
$$\chi(\mathcal{O}_X(D^{'}))\geq \chi(\mathcal{O}_X(2H-D))-\frac{1}{3}(H^2-4),$$
\noindent and hence, we have ${D^{'}}^2\geq \frac{1}{3}(H^2-4)>0$. This is a contradiction. Hence, we have ${D^{'}}^2>0$. By the same way as above, we have ${D^{'}}^2\geq (2H-D)^2=H^2-4$. By Hodge index theorem, we have
$$(H^2-3)^2<(H^2)(H^2-4)\leq (H.D^{'})^2\leq (H^2-1)^2.$$
\noindent Hence, we have $H.D^{'}=H^2-2$ or $H^2-1$. If $H.D^{'}=H^2-2$, then we have $(H-D^{'})^2\geq 0$ and $H.(H-D^{'})=2$. However, by Corollary 2.1 (iii), this contradicts to the assumption that $L$ is very ample. Assume that $H.D^{'}=H^2-1$. Note that, since $H.\Delta=1$, $\Delta$ is a $(-2)$-curve. If ${D^{'}}^2>H^2-4$, then we have $(H-D^{'})^2>-2$ and $H.(H-D^{'})=1$. By the same reason, this is a contradiction. Therefore, we have  ${D^{'}}^2=H^2-4$. Here, we note that, since $(H-D^{'})^2=-2$ and $H.(H-D^{'})=1$, the member of $|H-D^{'}|$ is a $(-2)$-curve. Since $(D^{'}+\Delta)^2=(2H-D)^2=H^2-4$, we have $D^{'}.\Delta=1$. Since $D^{'}.(2H-D)=D^{'}.(D^{'}+\Delta)=H^2-3$, we have $D^{'}.D=H^2+1$. Hence, we have $D.(H-D^{'})=-1<0$. This contradicts to the assumption that $|D|$ is base point free. Hence, $|2H-D|$ is base point free.

We consider the case where $|D|$ is not base point free. Let $\Delta$ be the base divisor of $|D|$, and let $D^{'}=D-\Delta$. Note that, since $H.D=H^2$, we have  $H.D^{'}\leq H^2-1$. Since $H^2\geq6$, by the same reason as above, we have ${D^{'}}^2>0$ and hence, we have $h^1(\mathcal{O}_X(D^{'}))=0$. Since $\chi(\mathcal{O}_X(D^{'}))\geq \chi(\mathcal{O}_X(D))$, we have ${D^{'}}^2\geq D^2=H^2-4$. Hence, by Hodge index theorem, we have
$$(H^2-3)^2<H^2(H^2-4)\leq (H.D^{'})^2\leq (H^2-1)^2.$$
\noindent Therefore, we have $H.D^{'}=H^2-2$ or $H^2-1$. If $H.D^{'}=H^2-2$, we have $(H-D^{'})^2\geq 0$ and $H.(H-D^{'})=2$. Since $L$ is very ample, by Corollary 2.1 (iii), this is a contradiction. If $H.D^{'}=H^2-1$, we have $H.\Delta=1$ and hence, $\Delta$ is a $(-2)$-curve. Hence, $D$ is a 1-connected divisor. In fact, since $D^2=H^2-4$, we have $2D^{'}.\Delta=H^2-2-{D^{'}}^2$. Since $D^{'}.\Delta\geq0$, we have ${D^{'}}^2=H^2-2$ or $H^2-4$. If ${D^{'}}^2=H^2-2$, then we have $(D^{'}-H)^2=0$. However, since $H.(H-D^{'})=1$, this contradicts to the ampleness of $L$. Hence, we have ${D^{'}}^2=H^2-4$ and hence, we have $D^{'}.\Delta=1$. Therefore, we have $h^1(\mathcal{O}_X(D))=0$. 

Next, we show that $h^1(\mathcal{O}_X(2H-D))=0$. Let $D^{'}$ and $\Delta$ be as above. Then we have $(H-\Delta)+(H-D^{'})=2H-D$ and $(H-\Delta).(H-D^{'})=1$. Moreover, since $H.(H-D^{'})=1$ and $(H-D^{'})^2=-2$, the member of $|H-D^{'}|$ is a $(-2)$-curve. Here, in order to show that $|2H-D|$ contains a 1-connected divisor, we show that $|H-\Delta|$ is base point free. We assume that it is not base point free and let $\Delta^{'}$ be the base divisor of it. Since $(H-\Delta)^2>0$, the divisor $H-\Delta-\Delta^{'}$ is not linearly equivalent to 0. Hence, we take a nonzero divisor $D^{''}\in |H-\Delta-\Delta^{'}|$. Note that, since $H.\Delta=1$, we have $H.D^{''}\leq H^2-2$. By the same reason as above, we have ${D^{''}}^2>0$ and hence, $h^1(\mathcal{O}_X(D^{''}))=0$. Therefore, we have ${D^{''}}^2\geq (H-\Delta)^2=H^2-4$. By Hodge index theorem, we have
$$(H^2-3)^2<H^2(H^2-4)\leq (H.D^{''})^2\leq (H^2-2)^2.$$
\noindent Hence, we have $H.D^{''}=H^2-2$. However, since $(H-D^{''})^2\geq0$ and $H.(H-D^{''})=2$, this contradicts to the very ampleness of $L$. Hence, $|H-\Delta|$ is base point free. Since $(H-\Delta)^2=H^2-4>0$, the general member of $|H-\Delta|$ is irreducible. Since $|2H-D|$ contains a 1-connected divisor, we have $h^1(\mathcal{O}_X(2H-D))=0$. Therefore, $\mathcal{O}_X(D)$ is ACM.

Next, we consider the case where $D^2=H^2-2$ and $H.D=H^2+1$. By Lemma 3.1, in order to show that $\mathcal{O}_X(D)$ is ACM, it is sufficient to show that 
$$h^1(\mathcal{O}_X(D))=h^1(\mathcal{O}_X(D-H))=h^1(\mathcal{O}_X(D-2H))=0.$$
First of all, we show that $|D|$ is base point free. Assume that it is not base point free, and let $\Delta$ be the base divisor of $D$ and let $D^{'}=D-\Delta$. Then we note that $H.D^{'}\leq H^2$. If ${D^{'}}^2=0$, then, by the same reason as above, we have a contradiction. Hence, we have ${D^{'}}^2>0$. This implies $\chi(\mathcal{O}_X(D^{'}))\geq\chi(\mathcal{O}_X(D))$. Hence, we have ${D^{'}}^2\geq D^2=H^2-2$. By Hodge index theorem, we have
$$(H^2-2)^2<H^2(H^2-2)\leq (H^2)({D^{'}}^2)\leq (H.D^{'})^2\leq (H^2)^2.$$
\noindent Hence, we have $H.D^{'}=H^2-1$ or $H^2$. If $H.D^{'}=H^2-1$, then $(H-D^{'})^2\geq 0$ and $H.(H-D^{'})=1$. This contradicts to Corollary 2.1 (iii). Hence, we have $H.D^{'}=H^2$. Note that, since $H.\Delta=1$, $\Delta$ is a $(-2)$-curve. Since $(H-D^{'})^2\geq -2$ and $H.(H-D^{'})=0$, we have $H\sim D^{'}$. However, since $D^2=H^2-2$ and $\Delta$ is a $(-2)$-curve, we have the contradiction $H.\Delta=0$. Therefore, $|D|$ is base point free and hence, we have $h^1(\mathcal{O}_X(D))=0$. Since $(H-D)^2=-4$, we have $\chi(\mathcal{O}_X(H-D))=0$. Since $|D|$ is base point free and $D.(D-H)=-3$, we have $|D-H|=\emptyset$. Hence, $\mathcal{O}_X(D)$ is initialized. Moreover, since $H.(H-D)=-1$, we have $|H-D|=\emptyset$. Hence, we have $h^1(\mathcal{O}_X(D-H))=0$.

We show that $h^1(\mathcal{O}_X(2H-D))=0$. If $|2H-D|$ is base point free, the assertion is satisfied. In fact, if $H^2=6$, then $(2H-D)^2=0$. Hence, by Proposition 2.2, there exists an elliptic curve $F$ and an integer $k\geq1$ such that $2H-D\sim kF$. In this case, since $H.(2H-D)=5$, we have $k=1$. Hence, we have $h^1(\mathcal{O}_X(2H-D))=0$. If $H^2\geq8$, we have $(2H-D)^2>0$, and hence, by the theorem of Bertini, we have the assertion. Therefore, we assume that it is not base point free. Let $\Delta$ be the base divisor of $|2H-D|$, and let $D^{'}\in|2H-D-\Delta|$. Note that, since $(2H-D)^2\geq0$, we have $D^{'}\neq0$.

First of all, we consider the case where $H^2=6$. Since $H.(2H-D)=5$, we have $H.D^{'}\leq 4$. If ${D^{'}}^2=0$, there exists an elliptic curve $F$ and an integer $k\geq1$ such that $D^{'}\sim kF$. Hence, by Corollary 2.1 (iii), we have $3k\leq H.D^{'}\leq 4$, and hence, we have $k=1$. Since $h^1(\mathcal{O}_X(D^{'}))=0$ and $(2H-D)^2={D^{'}}^2$, by Lemma 3.2, we have $h^1(\mathcal{O}_X(2H-D))=0$. If ${D^{'}}^2>0$, by Hodge index theorem, we have
$$(H.D^{'})^2\geq (H^2)({D^{'}})^2\geq 12.$$
\noindent We have $H.D^{'}=4$. Since we have $(H-D^{'})^2\geq0$ and $H.(H-D^{'})=2$, by Corollary 2.1 (iii), we have a contradiction. Hence, we have the assertion.

Next, we consider the case where $H^2=8$. Since $H.(2H-D)=7$, we have $H.D^{'}\leq 6$. If ${D^{'}}^2=0$, there exists an elliptic curve $F$ and an integer $k\geq1$ such that $D^{'}\sim kF$. We have $3k\leq H.D^{'}\leq 6$. Hence, we have $k=1$ or 2. If $k=1$, we have $h^1(\mathcal{O}_X(D^{'}))=0$ and hence, $\chi(\mathcal{O}_X(D^{'}))\geq \chi(\mathcal{O}_X(2H-D))$. This implies the contradiction ${D^{'}}^2\geq (2H-D)^2=2$. If $k=2$, we have $H.\Delta=1$ and hence, $\Delta$ is a $(-2)$-curve. Since $(D^{'}+\Delta)^2=(2H-D)^2=2$, we have $D^{'}.\Delta=2$. Hence, we have $h^1(\mathcal{O}_X(2H-D))=0$. If ${D^{'}}^2=2$, we have $h^1(\mathcal{O}_X(D^{'}))=0$ and ${D^{'}}^2=(2H-D)^2=2$. Hence, in this case, by Lemma 3.2, we have $h^1(\mathcal{O}_X(2H-D))=0$. Assume that ${D^{'}}^2\geq 4$. By Hodge index theorem, we have $32\leq(H^2)({D^{'}}^2)\leq (H.D^{'})^2$. Since $H.D^{'}=6$, we have $(H-D^{'})^2\geq 0$ and $H.(H-D^{'})=2$. By Corollary 2.1 (iii), this is a contradiction.

Finally, we consider the case where $H^2\geq 10$. We note that, since $H.(2H-D)=H^2-1$, we have $H.D^{'}\leq H^2-2$. Assume that ${D^{'}}^2=0$. Then there exists an elliptic curve $F$ and an integer $k\geq1$ such that $D^{'}\sim kF$. By Corollary 2.1 (iii), we have $3k\leq H.D^{'}\leq H^2-2$. Since
$$\chi(\mathcal{O}_X(D^{'}))\geq h^0(\mathcal{O}_X(D^{'}))-k+1\geq \chi(\mathcal{O}_X(2H-D))-\frac{1}{3}(H^2-5),$$
\noindent we have the contradiction ${D^{'}}^2\geq \frac{1}{3}(H^2-8)>0$. Since we have ${D^{'}}^2>0$, we have ${D^{'}}^2\geq (2H-D)^2=H^2-6$. By Hodge index theorem, we have 
$$(H^2-4)^2< H^2(H^2-6)\leq (H^2)({D^{'}}^2)\leq (H.D^{'})^2\leq (H^2-2)^2.$$
\noindent Hence, we have $H.D^{'}=H^2-3$ or $H^2-2$. 

Assume that $H.D^{'}=H^2-3$. If ${D^{'}}^2>H^2-6$, then $(H-D^{'})^2>0$. By Hodge index theorem, we have the contradiction 
$$9=(H.(H-D^{'}))^2\geq (H^2)(H-D^{'})^2\geq 20.$$
\noindent If ${D^{'}}^2=H^2-6$, we have $h^1(\mathcal{O}_X(D^{'}))=0$ and ${D^{'}}^2=(2H-D)^2$. Hence, by Lemma 3.2, we have $h^1(\mathcal{O}_X(2H-D))=0$. 

Assume that $H.D^{'}=H^2-2$. Since $H.\Delta=1$, $\Delta$ is a $(-2)$-curve. Since $(D^{'}+\Delta)^2=(2H-D)^2=H^2-6$, we have
$$2D^{'}.\Delta=H^2-4-{D^{'}}^2.$$
Since $D^{'}.\Delta\geq0$, we have ${D^{'}}^2 = H^2-4$ or $H^2-6$. If ${D^{'}}^2=H^2-4$, we have $(H-D^{'})^2=0$ and $H.(H-D^{'})=2$. However, by Corollary 2.1 (iii), this is a contradiction. Since ${D^{'}}^2=H^2-6$, we have $D^{'}.\Delta=1$. Therefore, $|2H-D|$ contains the 1-connected divisor $D^{'}+\Delta$. Hence, we have $h^1(\mathcal{O}_X(2H-D))=0$.

$\;$

We consider the case where $D^2\geq H^2$. Note that, by the proof of ${\rm{(i)}}\Rightarrow {\rm{(ii)}}$, we have $|H-D|=\emptyset$. Hence, by the assumption, we have
$$h^1(\mathcal{O}_X(D-H))=h^1(\mathcal{O}_X(D-2H))=0.$$
\noindent Since
$$\chi(\mathcal{O}_X(D-H))=0\text{ and }\chi(\mathcal{O}_X(2H-D))\geq0,$$
\noindent we have $D^2\leq 2H^2-4$. Hence, we have 
$$H.D=\frac{1}{2}(D^2+H^2+4)\leq\frac{3}{2}H^2.$$
\noindent Therefore, by Lemma 3.1, it is sufficient to show that $h^1(\mathcal{O}_X(D))=0$. Since $D^2>0$, we show that $|D|$ is base point free. Assume that $|D|$ is not base point free, and let $\Delta$ be the base divisor of $|D|$. Let $D^{'}=D-\Delta$. Since $D^2\leq 2H^2-4$ and $L$ is ample, we have
$$1\leq H.D^{'}<H.D\leq\frac{3}{2}H^2.$$
\noindent Assume that ${D^{'}}^2=0$. Then there exists an elliptic curve $F$ and an integer $k$ such that $D^{'}\sim kF$. Since $L$ is very ample, by Corollary 2.1 (iii), we have
$$H.D^{'}=kH.F\geq 3k.$$
\noindent Since we have $k\leq\frac{1}{2}H^2-1$, we have
$$\chi(\mathcal{O}_X(D^{'}))=h^0(\mathcal{O}_X(D^{'}))-k+1\geq \chi(\mathcal{O}_X(D))-\frac{1}{2}H^2+2.$$
\noindent Hence, we have
$${D^{'}}^2\geq D^2-H^2+4\geq 4.$$
\noindent This is a contradiction. Hence, we have ${D^{'}}^2>0$. Since $h^1(\mathcal{O}_X(D^{'}))=0$, we have
$${D^{'}}^2\geq D^2\geq H^2.$$ 
\noindent By Hodge index theorem, we have
$$(H^2)^2\leq (D^2)(H^2)\leq ({D^{'}}^2)(H^2)\leq (H.D^{'})^2.$$
\noindent Hence, we have $H.(D^{'}-H)\geq 0$. Since $H.\Delta>0$, by the assumption, we have 
$$(D^{'}-H)^2>(D-H)^2=-4.$$
\noindent This implies $|D^{'}-H|\neq\emptyset$. However, this contradicts to the assumption that $|D-H|=\emptyset$. Hence, $|D|$ is base point free. The assertion holds. $\hfill\square$

\section{Proof of Theorem 1.1}

In this section, we give a proof of Theorem 1.1, by using Theorem 3.2.

\smallskip

\smallskip

{\it{Proof of Theorem 1.1}}. First of all, let $H^{'}\in |\pi^{\ast}\mathcal{O}_{\mathbb{P}^2}(1)|$, and we note that $H\sim 3H^{'}$ and hence, $3|H.D$.

\smallskip

\smallskip

\noindent ${\rm{(i)}}\Rightarrow {\rm{(ii)}}$ We note that, since $\mathcal{O}_X(D)$ is initialized and $h^1(\mathcal{O}_X(D))=0$, we have $D^2\geq-2$. Assume that $D^2=-2$. If $|H-D|\neq\emptyset$, by the assumption that $\mathcal{O}_X(D)$ is ACM, we have $\chi(\mathcal{O}_X(H-D))\geq1$ and hence, we have $H.D\leq9$. Since $3|H.D$, we have $H.D=3,6$ or $9$. Since $h^1(\mathcal{O}_X(H-D))=0$ and $|D-H|=\emptyset$, if $|H-D|=\emptyset$, we have $\chi(\mathcal{O}_X(H-D))=0$. Hence, we have $H.D=10$. This contradicts to $3|H.D$.

Assume that $D^2=0$. Then, by the ampleness of $H^{'}$ and Corollary 2.1 (i), we have $H^{'}.D\geq 3$. Hence, we have $H.D\geq9$. Assume that $|H-D|\neq\emptyset$. Since $h^1(\mathcal{O}_X(H-D))=0$, we have $\chi(\mathcal{O}_X(H-D))\geq1$. Hence, we have $H.D\leq 10$. Since $3|H.D$, we have $H.D=9$. If $|H-D|=\emptyset$, by the same reason as above, we have a contradiction. 

Assume that $D^2=2$. By Corollary 2.1 (i), we have $D\sim H^{'}$ or $H.D\geq9$. If $D\sim H^{'}$, then we have $H.D=6$. We assume that the latter case occurs. If $|H-D|\neq\emptyset$, then we have $\chi(\mathcal{O}_X(H-D))\geq 1$. In this case, since $H.D\leq 11$, we have $H.D=9$. If $|H-D|=\emptyset$, then we have $\chi(\mathcal{O}_X(H-D))=0$. Hence, we have $H.D=12$. Conversely, if $H.D=12$, then we have $|H-D|=\emptyset$, by the assumption.

Assume that $D^2=4$. By Corollary 2.1 (i), we have $H.D\geq9$. If $|H-D|\neq\emptyset$, by the same reason as above, we have $H.D\leq 12$. In this case, we have $H.D=9$ or 12. If $|H-D|=\emptyset$, by the same reason as above, we have a contradiction.

Assume that $D^2=6$. We have $H.D\geq9$. If $|H-D|\neq\emptyset$, by the same reason as above, we have $H.D\leq 13$. Hence, we have $H.D=9$ or 12. If $H.D=9$, we have $(H-D)^2=6$ and $H.(H-D)=9$. However, by Hodge index theorem, we have the contradiction
$$81=(H.(H-D))^2\geq (H^2)(H-D)^2=108.$$
\noindent If $H.D=12$, we have $(H-D)^2=0$ and $H^{'}.(H-D)=2$. This contradicts to Corollary 2.1 (i). If $|H-D|=\emptyset$, by the same reason as above, we have the contradiction $H.D=14$.

Assume that $D^2=8$. By Hodge index theorem, we have 
$$(H.D)^2\geq (H^2)(D^2)=144.$$ 
\noindent If $|H-D|\neq\emptyset$, by the same reason as above, we have $H.D\leq 14$. Hence, we have $H.D=12$. If $|H-D|=\emptyset$, we have $H.D=15$.

Assume that $D^2=10$. By Hodge index theorem, we have 
$$(H.D)^2\geq (H^2)(D^2)=180.$$
\noindent If $|H-D|\neq\emptyset$, we have $H.D\leq 15$. Hence, we have $H.D=15$. If $|H-D|=\emptyset$, we have the contradiction $H.D=16$.

Assume that $D^2=12$. By the same reason as above, we have $H.D>14$. If $|H-D|\neq \emptyset$, we have $H.D\leq16$ and hence, we have $H.D=15$. In this case, we have $(H-D)^2=0$ and $H^{'}.(H-D)=1$. This contradicts to Corollary 2.1 (i). If $|H-D|=\emptyset$, by the same reason as above, we have the contradiction $H.D=17$.

We consider the case where $D^2\geq 14$. Assume that $D^2=14$ or 16. Since $3|H.D$, by Theorem 3.2, we have $D^2=14$ and $H.D=18$. Assume that $D^2\geq 18$. We note that, by the proof of Theorem 3.2 ${\rm{(i)}}\Rightarrow {\rm{(ii)}}$, we have $|H-D|=\emptyset$. Since $\mathcal{O}_X(D)$ is ACM and initialized, we have $\chi(\mathcal{O}_X(H-D))=0$ and $\chi(\mathcal{O}_X(2H-D))\geq0$. Hence, we have $D^2\leq 32$. Since $3|H.D$, by Theorem 3.2, we have $D^2=20,\;26$ or 32.

\smallskip

\smallskip

(ii) $\Rightarrow$ (i) If $D^2\geq 14$, the assertion follows from Theorem 3.2. Hence, we consider the case where $-2\leq D^2\leq 10$. Since $H.(D-H)<0$, by the ampleness of $L$, we have $|D-H|=\emptyset$ and hence, $\mathcal{O}_X(D)$ is initialized. Hence, by Lemma 3.1, it is sufficient to show that $h^1(\mathcal{O}_X(D))=h^1(\mathcal{O}_X(H-D))=0$.

We consider the case where $D^2=-2$. Assume that $H.D=3$. Since $H^{'}.D=1$, $D$ is irreducible. Hence, we have $h^1(\mathcal{O}_X(D))=0$. Since $(H-D)^2=10$ and $H.(H-D)=15$, we have $|H-D|\neq\emptyset$. In order to show that $h^1(\mathcal{O}_X(H-D))=0$, we show that $|H-D|$ is base point free. Assume that it is not base point free and let $\Delta$ be the base divisor of $|H-D|$. Since $(H-D)^2>0$, we take a nonzero divisor $D^{'}\in|H-D-\Delta|$. Then we note that $H.D^{'}\leq 14$. Assume that ${D^{'}}^2=0$. Then, by Proposition 2.2, there exists an elliptic curve $F$ and an integer $k\geq 1$ such that $D^{'}\sim kF$. By Corollary 2.1 (i), we have $H^{'}.F\geq3$ and hence, we have $9k\leq H.D^{'}\leq 14$. Since we have $k=1$, we have $h^1(\mathcal{O}_X(D^{'}))=0$. Hence, by comparing $\chi(\mathcal{O}_X(H-D))$ and $\chi(\mathcal{O}_X(D^{'}))$, we have the contradiction ${D^{'}}^2\geq (H-D)^2=10$. Since ${D^{'}}^2>0$, we have ${D^{'}}^2\geq (H-D)^2=10.$ Hence, by Hodge index theorem, we have
$$180\leq(H^2)({D^{'}}^2)\leq (H.D^{'})^2\leq 196.$$
\noindent Hence, we have $H.D^{'}=14$. However, this contradicts to $3|H.D^{'}$. Therefore, $|H-D|$ is base point free and hence, we have $h^1(\mathcal{O}_X(H-D))=0$.

Assume that $H.D=6$. Since $H^{'}.D=2$, by Corollary 2.1 (i), the movable part of $|D|$ is empty. Hence, we have $h^0(\mathcal{O}_X(D))=1$ and $h^1(\mathcal{O}_X(D))=0$. Since we have $(H-D)^2=4$ and $H.(H-D)=12$, if $|H-D|$ is base point free, then, by the theorem of Bertini, we have $h^1(\mathcal{O}_X(H-D))=0$. Hence, we assume that it is not base point free. Let $\Delta$ be the base divisor of it and take $D^{'}\in|H-D-\Delta|$. If ${D^{'}}^2=0$, by the same reason as above, there exists an elliptic curve $F$ such that $D^{'}\sim F$. Hence, by comparing $\chi(\mathcal{O}_X(H-D))$ and $\chi(\mathcal{O}_X(D^{'}))$, we have the contradiction ${D^{'}}^2\geq 4$. Hence, by the same reason as above, we have ${D^{'}}^2\geq4$. By Hodge index theorem, we have
$$72\leq(H^2)({D^{'}}^2)\leq (H.D^{'})^2\leq 121.$$
\noindent Hence, we have $H.D^{'}=9$ and ${D^{'}}^2=4$. Since we have $H^{'}.\Delta=1$, $\Delta$ is a $(-2)$-curve. Moreover, since 
$$(D^{'}+\Delta)^2=(H-D)^2=4,$$
\noindent we have $D^{'}.\Delta=1$. Since $|H-D|$ contains a 1-connected divisor $D^{'}+\Delta$, we have $h^1(\mathcal{O}_X(H-D))=0$.

Assume that $H.D=9$. Since $D^2=(H-D)^2=-2$ and $H.D=H.(H-D)=9$, by symmetry, it is sufficient to show that $h^1(\mathcal{O}_X(D))=0$. Since $D^2=-2$, the base divisor $\Delta$ of $|D|$ is not zero. Let $D^{'}=D-\Delta$. If $D^{'}\neq 0$, then, by Corollary 2.1 (i), we have $H^{'}.D^{'}\geq3$. In fact, if $H^{'}\sim D^{'}$, then we have $H^{'}.\Delta=1$ and hence, $\Delta$ is a $(-2)$-curve. However, since $D^2=-2$, we have the contradiction $H^{'}.\Delta=-1$. Hence, we have $H.D^{'}=9$ and $H.\Delta=0$. This contradicts to the ampleness of $L$. Hence, $D^{'}=0$ and hence, we have $h^0(\mathcal{O}_X(D))=1$. Since $D^2=-2$, we have
$$h^1(\mathcal{O}_X(D))=-\chi(\mathcal{O}_X(D))+1=0.$$

We consider the case where $D^2=0$. Since $H.(H-D)=H.D=9$ and $(H-D)^2=D^2=0$, we only show that $h^1(\mathcal{O}_X(D))=0$. Since $H^{'}.D=3$, by Corollary 2.1 (ii), $|D|$ is base point free and $D$ is irreducible. Hence, we have $h^1(\mathcal{O}_X(D))=0$.

We consider the case where $D^2=2$. If $H.D=6$, since we have $(H^{'}-D)^2=0$ and $H^{'}.(H^{'}-D)=0$, we have $D\sim H^{'}$. Hence, it is clear that $h^1(\mathcal{O}_X(D))=h^1(\mathcal{O}_X(H-D))=0$. If $H.D=9$, we have $(H-D)^2=D^2=2$ and $H.(H-D)=H.D=9$. Hence, we only show that $h^1(\mathcal{O}_X(D))=0$. If $D$ is base point free, then the assertion follows immediately. Hence, we assume that it is not base point free. By the proof of Corollary 2.1 (ii), there exists  a $(-2)$-curve $\Gamma$ such that $D\sim H^{'}+\Gamma$ and $H^{'}.\Gamma=1$. Since $D$ is 1-connected, we have $h^1(\mathcal{O}_X(D))=0$.

Assume that $H.D=12$. First of all, since $(H-D)^2=-4$ and $|D-H|=|H-D|=\emptyset$, we have $h^1(\mathcal{O}_X(H-D))=0$. In order to show that $h^1(\mathcal{O}_X(D))=0$, we show that $|D|$ is base point free. Assume that it is not base point free. Let $\Delta$ be the base divisor of it, and let $D^{'}=D-\Delta$. Then, by Corollary 2.1 (i), we have $D^{'}\sim H^{'}$ or $H^{'}.D^{'}\geq 3$. However, the first case does not occur. In fact, if $D^{'}\sim H^{'}$, then we have $H^{'}.\Delta=2$. Since $D^2=2$, we have $\Delta^2=-4$. Hence, there exist $(-2)$-curves $\Gamma_1$ and $\Gamma_2$ such that $\Delta=\Gamma_1+\Gamma_2$, $\Gamma_1.\Gamma_2=0$, and $H^{'}.\Gamma_i=1\;(i=1,2)$. Hence, we have $H-D\sim 2H^{'}-\Gamma_1-\Gamma_2$. Since $(H^{'}-\Gamma_i)^2=-2$ and $H^{'}.(H^{'}-\Gamma_i)=1$, the member of $|H^{'}-\Gamma_i|$ is a $(-2)$-curve. This contradicts to the assumption that $|H-D|=\emptyset$. Hence, we have $H^{'}.D^{'}\geq 3$. Since we have $H.D^{'}=9$ and $H^{'}.\Delta=1$, $\Delta$ is a $(-2)$-curve. If ${D^{'}}^2=0$, there exists an elliptic curve $F$ such that $D^{'}\sim F$. Hence, by comparing $\chi(\mathcal{O}_X(D^{'}))$ and $\chi(\mathcal{O}_X(D))$, we have the contradiction ${D^{'}}^2\geq D^2=2$. Hence, we have ${D^{'}}^2>0$. If ${D^{'}}^2=2$, since $(D^{'}-H^{'})^2=-2$ and $H^{'}.(D^{'}-H^{'})=1>0$, we have $|D^{'}-H^{'}|\neq\emptyset$. However, we have the contradiction $D^{'}.(D^{'}-H^{'})=-1<0$. If ${D^{'}}^2=4$, since $(D^{'}-H^{'})^2=0$ and $H^{'}.(D^{'}-H^{'})=1$, we have a contradiction, by Corollary 2.1 (i). If ${D^{'}}^2\geq6$, by Hodge index theorem, we have
$$108\leq({D^{'}}^2)(H^2)\leq (H.D^{'})^2\leq 121.$$
\noindent However, this contradicts to $3|H.D^{'}$. Hence, $|D|$ is base point free and hence, $h^1(\mathcal{O}_X(D))=0$.

We consider the case where $D^2=4$. Assume that $H.D=9$. Since $(H-D)^2=D^2=4$ and $H.(H-D)=H.D=9$, it is sufficient to show that $h^1(\mathcal{O}_X(D))=0$. Since $H^{'}.D=3$, by Corollary 2.1 (ii), $|D|$ is base point free. Hence, by the theorem of Bertini, we have $h^1(\mathcal{O}_X(D))=0$. If $H.D=12$, we have $(H-D)^2=-2$ and $H.(H-D)=6$. In this case, by the case of (a) in Theorem 1.1 (we have already proved it), $\mathcal{O}_X(H-D)$ is ACM. Hence, $\mathcal{O}_X(D)$ is also ACM.

We consider the case where $D^2=8$. If $H.D=12$, then we have $(H-D)^2=2$ and $H.(H-D)=6$. In this case, by the case of (c) in Theorem 1.1 (we have already proved it), $\mathcal{O}_X(H-D)$ is ACM. Hence, $\mathcal{O}_X(D)$ is also ACM. 

Assume that $H.D=15$. First of all, in order to show that $h^1(\mathcal{O}_X(D))=0$, we show that $|D|$ is base point free. Assume that $|D|$ is not base point free and let $\Delta$ be the base divisor of $|D|$. Let $D^{'}\in |D-\Delta|$. If ${D^{'}}^2=0$, by the same argument as above, we have a contradiction. Hence, we have ${D^{'}}^2>0$ and hence, we have $h^1(\mathcal{O}_X(D^{'}))=0$. By comparing $\chi(\mathcal{O}_X(D))$ and $\chi(\mathcal{O}_X(D^{'}))$, we have ${D^{'}}^2\geq D^2=8$. Since $H.D=15$, we note that $H.D^{'}\leq14$. By Hodge index theorem, we have 
$$144\leq ({D^{'}}^2)(H^2)\leq (H.D^{'})^2\leq 196.$$
\noindent Hence, in this case, we have $H.D^{'}=12$ and $H^{'}.\Delta=1$. Therefore, $\Delta$ is a $(-2)$-curve. If ${D^{'}}^2=8$, since $(D^{'}+\Delta)^2=D^2=8$, we have $D^{'}.\Delta=1$. In this case, we also have $(D^{'}-H^{'})^2=H^{'}.(D^{'}-H^{'})=2$. Hence, by Corollary 2.1 (i), we have $D^{'}\sim2H^{'}$. However, this contradicts to $D^{'}.\Delta=1$. If ${D^{'}}^2\geq10$, by Hodge index theorem, we have
$$180\leq(H^2)({D^{'}}^2)\leq (H.D^{'})^2\leq 196.$$
\noindent This contradicts to $3|H.D^{'}$. Therefore, $|D|$ is base point free and hence, we have $h^1(\mathcal{O}_X(D))=0$. Next, we show that $h^1(\mathcal{O}_X(H-D))=0$. Since $(H-D)^2=-4$ and $H^{'}.(H-D)=1$, we have $|H-D|=\emptyset$. Since we already have $|D-H|=\emptyset$, we have $$h^1(\mathcal{O}_X(H-D))=-\chi(\mathcal{O}_X(H-D))=0.$$

Assume that $D^2=10$. Since $(H-D)^2=-2$ and $H.(H-D)=3$, by the case of (a) in Theorem 1.1, $\mathcal{O}_X(H-D)$ is ACM. Hence, $\mathcal{O}_X(D)$ is also ACM. $\hfill\square$

\pagebreak

\section{K3 surfaces and 2-elementary lattices}

In this section, we recall the definition of a 2-elementary lattice and some basic results about a K3 surface whose N$\rm{\acute{e}}$ron-Severi lattice is a 2-elementary lattice.

\begin{df} A lattice $S$ is called a {\it{2-elementary}} lattice if there exists a non-negative integer $a$ such that $S^{\ast}/S\cong (\mathbb{Z}/2\mathbb{Z})^a$, where $S^{\ast}:=\Hom(S,\mathbb{Z})$. \end{df}

\begin{df} For a 2-elementary lattice $S$, we define

\[\delta_{S}=\left\{
\begin{array}{ll}
\displaystyle{0}&
\qquad x^2\in\mathbb{Z}, \forall x\in S^{\ast}\\[0.2cm]
\displaystyle{1}&\qquad \text{otherwise.}
\end{array}\right. \]

\end{df}

We note that, by the classification of 2-elementary lattices (cf. [Ni], Theorem 4.3.2), we have the following assertion.

\begin{prop}Let $S$ be a hyperbolic, even, 2-elementary lattice, and let $a$ and $\delta$ be as above. We assume that the rank of $S$ is $a$ and $1\leq a\leq 9$. Then, if $\delta_S=0$, then $a=2$ and $S=U(2)$. On the other hand, if $\delta_S=1$, then $S =< 2 >\oplus A_1^{\oplus a-1}$. \end{prop}

\noindent {\bf{Remark 5.1}}. By Proposition 5.1, if $X$ is a K3 surface whose N$\rm{\acute{e}}$ron-Severi lattice is a 2-elementary lattice with $\rho(X)=a$ and $1\leq a\leq9$, then the intersection number of any two divisors on $X$ is even.

\smallskip

\smallskip

Let $X$ be a K3 surface, and let $\omega_X$ be a nowhere vanishing holomorphic 2-form on $X$. Then we call an automorphism $\varphi$ {\it{non-symplectic}} if the action of $\varphi$ on the complex vector space $H^0(K_X)$ is not trivial. If the N$\rm{\acute{e}}$ron-Severi lattice $S_X$ of $X$ is a 2-elementary lattice, then there exists a unique non-symplectic involution $\theta$ which acts trivially on $S_X$. We call it the canonical involution. It is well known that the fixed locus of $\theta$ forms a divisor on $X$, and the concrete description of it is given as follows (cf. [Ni], Theorem 4.2.2).

\smallskip

\smallskip

\begin{thm} Let $X$ be a K3 surface whose N$\acute{e}$ron-Severi lattice is a 2-elementary lattice, and let $\theta$ be as above. Then the set of fixed points $X^{\theta}$ has the form

\[X^{\theta}=\left\{
\begin{array}{ll}
\displaystyle{\phi},&
 (\rho(X), a, \delta_{S_X})=(10,10,0)\\[0.2cm]
\displaystyle{C_1^{(1)}+C_2^{(1)},}&
 (\rho(X),a,\delta_{S_X})=(10,8,0) \\[0.2cm]
\displaystyle{C^{(g)}+\sum_{1\leq i\leq k}E_i,}&\qquad \rm{otherwise.}
\end{array}\right. \]

\end{thm}

\noindent Here $a$ is the minimal number of generators of ${S_X}^{\ast}/S_X$ and $k=(\rho(X)-a)/2$. We denote by $C^{(g)}$ a curve of genus $g$, where $g=(22-\rho(X)-a)/2$, and by $E_i$ a smooth rational curve. We note that $C^{(g)}$ and $E_i\;(1\leq i\leq k)$ do not intersect each other. If $g\geq2$, then we say that the involution $\theta$ on $X$ is of {\it{elliptic type}}.

\smallskip

\smallskip

Let $X$ be a K3 surface as in Theorem 5.1 and $\theta$ be the canonical involution on $X$. Then the quotient surface $Y:=X/\langle \theta\rangle$ is a smooth surface with $\rho(X)=\rho(Y)$ (if $X^{\theta}=\emptyset$, then $Y$ is an Enriques surface, otherwise $Y$ is a rational surface).  $X^{\theta}$ is ample as a divisor on $X$ if and only if $X^{\theta}$ is irreducible and $\theta$ is of elliptic type, that is,  $\rho(X)=a$ and $1\leq a\leq9$. In particular, by Theorem 5.1, if $\rho(X)=a=9$, the set of fixed points $X^{\theta}$ of the canonical involution $\theta$ is a smooth genus 2 curve. Hence, the linear system $|X^{\theta}|$ gives a double cover $\pi:X\rightarrow \mathbb{P}^2$ branched along a smooth sextic, and $\pi^{\ast}\mathcal{O}_{\mathbb{P}^2}(1)$ is an ample line bundle of sectional genus 2. Hence, we have the following result.

\begin{prop} Let $X$ be a K3 surface whose N$\acute{e}$ron-Severi lattice is a 2-elementary lattice with $\rho(X)=a=9$. Let $D$ be a nonzero effective divisor on $X$, and let $L=\mathcal{O}_X(3X^{\theta})$. Then the following conditions are equivalent.

\smallskip

\smallskip

\noindent {\rm{(i)}} $\mathcal{O}_X(D)$ is ACM and initialized with respect to $L$.

\smallskip

\noindent {\rm{(ii)}} For $H\in |L|$, one of the following cases occurs.

{\rm{(a)}} $D^2=-2$ and $H.D=6$.

{\rm{(b)}} $D^2=2$ and $H.D=6$.

{\rm{(c)}} $D^2=4$ and $H.D=12$.

{\rm{(d)}} $D^2=8$ and $H.D=12$.

{\rm{(e)}} $D^2=14$ and $H.D=18$.

\smallskip

\noindent {\rm{(iii)}} One of the following cases occurs.

{\rm{(f)}} There exists a $(-2)$-curve $\Gamma$ such that $D\sim \Gamma,\;3X^{\theta}-\Gamma$ or $4X^{\theta}-\Gamma$.

{\rm{(g)}} $D\sim rX^{\theta}\;(r=1,2)$. \end{prop}

{\it{Proof}}. (i)$\Rightarrow$ (ii) Note that, by Remark 5.1, $H.D$ is even. Hence, it is sufficient to show that the case where $D^2=2$ and $H.D=12$ and the case where $D^2=26$ as in Theorem 1.1 do not occur. 

Assume that $D^2=2$, $H.D=12$ and $|H-D|=\emptyset$. By the proof of Theorem 1.1, $|D|$ is base point free. Hence, we can assume that $D$ is a smooth curve. Since $X^{\theta}.D\neq2$, we have $D\notin |X^{\theta}|$. Since $\theta$ acts trivially on $S_X$, if $\theta(D)\neq D$, we have $X^{\theta}.D\leq \theta(D).D=D^2=2$, otherwise, since $\theta$ acts on $D$, by the Hurwitz's formula, we have
$$2=D^2=4(g(D/\langle\theta\rangle)-1)+X^{\theta}.D.$$
\noindent In the latter case, since $g(D/\langle\theta\rangle)\leq1$, we have $X^{\theta}.D=2$ or 6. However, this contradicts to $H.D=12$. 

Assume that $D^2=26$. By Theorem 1.1, we have $H.D=24$. Note that, since $\mathcal{O}_X(D)$ is initialized, $|D-H|=\emptyset$. Since $(2H-D)^2=2$ and $X^{\theta}.(2H-D)=4$, by the same reason as above, $|2H-D|$ is base point free. Hence, the general member of it is a smooth genus 2 curve. By the same reason as above, we have the contradiction $X^{\theta}.(2H-D)=2$ or 6. 

\smallskip

\smallskip

(ii)$\Rightarrow$ (iii) Assume that $D^2=-2$ and $H.D=6$. Since $X^{\theta}.D=2$, by the ampleness of $X^{\theta}$ and Remark 5.1, $D$ is a $(-2)$-curve. 

Assume that $D^2=2$ and $H.D=6$. Since we have $X^{\theta}.D=2$, by the proof of Corollary 2.1 (i), we have $D\sim X^{\theta}$. 

Assume that $D^2=4$ and $H.D=12$. Since $X^{\theta}.D=4$, we have $(3X^{\theta}-D)^2=-2$ and $X^{\theta}.(3X^{\theta}-D)=2$. Hence, the member of $|3X^{\theta}-D|$ is a $(-2)$-curve. Hence, for $\Gamma\in |3X^{\theta}-D|$, we have $D\sim 3X^{\theta}-\Gamma$. 

If $D^2=8$ and $H.D=12$, since $(H-D)^2=2$ and $H.(H-D)=6$, we have $H-D\sim X^{\theta}$. Hence, in this case, we have $D\sim 2X^{\theta}$.

 Assume that $D^2=14$ and $H.D=18$. Since $(4X^{\theta}-D)^2=-2$ and $X^{\theta}.(4X^{\theta}-D)=2$, by the same reason as above, the member of $|4X^{\theta}-D|$ is a $(-2)$-curve. Hence, for $\Gamma\in |4X^{\theta}-D|$, we have $D\sim 4X^{\theta}-\Gamma$.
 
\smallskip

\smallskip

(iii)$\Rightarrow$(i) By Theorem 1.1, the assertion follows immediately. $\hfill\square$

\section{Extensions of vector bundles}

\noindent First of all, in order to consider the construction of families of semistable indecomposable ACM bundles of higher rank, we recall the criterion for semistability of vector bundles and some results about extensions of vector bundles.

\begin{df} Let $X$ be a smooth projective variety. Then a vector bundle $\mathcal{E}$ on $X$ is called {\it{semistable}} if, for any nonzero coherent subsheaf $\mathcal{F}$ of $\mathcal{E}$, there exists an integer $n_0\geq1$ such that, if $n\geq n_0$, then the inequality
$$\chi(\mathcal{F}(n))/\rk(\mathcal{F})\leq \chi(\mathcal{E}(n))/\rk(\mathcal{E})$$
\noindent holds, where $\rk(\mathcal{G})$ is the rank of $\mathcal{G}$.\end{df}

\noindent In Definition 6.1, a vector bundle $\mathcal{E}$ is said to be {\it{stable}} if, for any nonzero subsheaf $\mathcal{F}$ and any integer $n>>0$, the equality of the inequality does not hold.  We can easily see that, by Definition 6.1, any line bundle is semistable. In this paper, we construct semistable bundles of higher rank, by using extensions of line bundles. Hence, we recall the following result about the semistability of a vector bundle.

\begin{lem} {\rm{(cf. [Mar], Lemma 1.4)}} Let $X$ be as in Definition 6.1, and let
$$ 0\longrightarrow\mathcal{E}^{'}\longrightarrow\mathcal{E}\longrightarrow\mathcal{E}^{''}\longrightarrow0$$
\noindent be an exact sequence of vector bundles on $X$ such that, for any $n\in\mathbb{N}$,
$$\chi(\mathcal{E}(n))/\rk(\mathcal{E}) = \chi(\mathcal{E}^{'}(n))/\rk(\mathcal{E}^{'}) = \chi(\mathcal{E}^{''}(n))/\rk(\mathcal{E}^{''}).$$
\noindent Then $\mathcal{E}$ is semistable if and only if $\mathcal{E}^{'}$ and $\mathcal{E}^{''}$ are semistable.
\end{lem}

\begin{df} Let $X$ be as in Definition 6.1. Let $\mathcal{F}$ and $\mathcal{G}$ be vector bundles on $X$, and let $\mathcal{E}$ and $\mathcal{E}^{'}$ be extensions of $\mathcal{G}$ by $\mathcal{F}$. Then we say $\mathcal{E}$ and $\mathcal{E}^{'}$ are {\rm{weak equivalent}} if there exist isomorphisms of vector bundles $\mathcal{F}\longrightarrow \mathcal{F}$, $\mathcal{G}\longrightarrow \mathcal{G}$ and $\psi:\mathcal{E}\longrightarrow\mathcal{E}^{'}$ such that the following diagram commutes:

\begin{equation}
\begin{CD}
0 @>>> \mathcal{F} @>>> \mathcal{E} @>>>\mathcal{G} @>>>0\\
@.             @VVV           @VV\text{$\psi$}V           @VVV\\
0 @>>> \mathcal{F} @>>> \mathcal{E}^{'} @>>>\mathcal{G} @>>>0.
\end{CD}
\end{equation} 

\noindent If $\mathcal{E}$ and $\mathcal{E}^{'}$ are weak equivalent, then we will write $\mathcal{E}\sim_w\mathcal{E}^{'}$.   

\end{df}

\noindent In Definition 6.2, we say the two extensions $\mathcal{E}$ and $\mathcal{E}^{'}$ are {\it{equivalent}} if we can take the isomorphism $\psi:\mathcal{E}\longrightarrow\mathcal{E}^{'}$ such that the two isomorphisms $\mathcal{F}\longrightarrow \mathcal{F}$ and $\mathcal{G}\longrightarrow \mathcal{G}$ as in the above diagram are identity. It is well known that equivalent classes of extensions of $\mathcal{G}$ by $\mathcal{F}$ correspond bijectively to the elements of $\Ext^1(\mathcal{G},\mathcal{F})$. Hence, the weak equivalence of two extensions $\mathcal{E}$ and $\mathcal{E}^{'}$ of $\mathcal{G}$ by $\mathcal{F}$ induces an equivalent relation of the two elements $[\mathcal{E}]$ and $[\mathcal{E}^{'}]$ of $\Ext^1(\mathcal{G},\mathcal{F})$ corresponding to $\mathcal{E}$ and $\mathcal{E}^{'}$ respectively. We also denote it by $[\mathcal{E}]\sim_w[\mathcal{E}^{'}]$.

\begin{df} For a given smooth variety $X$, a vector bundle $\mathcal{E}$ is called simple if $\Hom(\mathcal{E},\mathcal{E})\cong\mathbb{C}$.\end{df}

\noindent We note that simple vector bundles are indecomposable. First of all, we recall the following results about extensions of simple vector bundles. 

\begin{prop} {\rm{(cf. [P-T], Proposition 5.1.3.)}} Let $X$ be a smooth projective variety. Let $\mathcal{F}_1,\cdots,\mathcal{F}_{r+1}\;(r\geq1)$ be simple vector bundles such that
$$\Hom(\mathcal{F}_i,\mathcal{F}_j)=0\;(i\neq j).$$
\noindent Moreover, let $\mathcal{F}:=\bigoplus_{1\leq i\leq r}\mathcal{F}_i$ and let
$$ U=\Ext^1(\mathcal{F}_{r+1},\mathcal{F}_1)\backslash \{0\}\times\cdots\times\Ext^1(\mathcal{F}_{r+1},\mathcal{F}_r)\backslash\{0\}\subset\Ext^1(\mathcal{F}_{r+1},\mathcal{F}).$$
\noindent Then we have the following results.

\smallskip

\smallskip

\noindent {\rm{(i)}} If an extension $\mathcal{E}$ of $\mathcal{F}_{r+1}$ by $\mathcal{F}$ is simple, then $[\mathcal{E}]\in U$. 

\noindent {\rm{(ii)}} For classes $[\mathcal{E}],\;[\mathcal{E}^{'}]\in U$ of extensions $\mathcal{E}$ and $\mathcal{E}^{'}$ of $\mathcal{F}_{r+1}$ by $\mathcal{F}$, 
$$\Hom(\mathcal{E},\mathcal{E}^{'})\neq0\Longleftrightarrow [\mathcal{E}]\sim_w[\mathcal{E}^{'}].$$
\end{prop}

\smallskip

\smallskip

\noindent  Since $\Ext^1(\mathcal{F}_{r+1},\mathcal{F})\cong \bigoplus_{1\leq i\leq r}\Ext^1(\mathcal{F}_{r+1},\mathcal{F}_i)$, for two elements $[\mathcal{E}],\;[\mathcal{E}^{'}]\in U$ as in Proposition 6.1 (ii), we can denote them as
$$[\mathcal{E}]=(\eta_1,\cdots,\eta_r)\text{ and }[\mathcal{E}^{'}]=(\xi_1,\cdots,\xi_r)\;(\eta_i,\xi_i\in\Ext^1(\mathcal{F}_{r+1},\mathcal{F}_i)\backslash\{0\}).$$
\noindent In Proposition 6.1, Joan Pons-Llopis and Fabio Tonini [P-T] have also proved that the condition as in (ii) is equivalent to the condition that there exists $\omega_i\in\mathbb{C}\backslash\{0\}$ such that $\eta_i=\omega_i\xi_i$, for any $1\leq i\leq r$. Therefore, we have
$$U/\sim_w\;\cong\;\mathbb{P}(\Ext^1(\mathcal{F}_{r+1},\mathcal{F}_1))\times\cdots\times\mathbb{P}(\Ext^1(\mathcal{F}_{r+1},\mathcal{F}_r)).$$

\section{Families of ACM bundles of higher rank}

\noindent In this section, we deal with the K3 surface of genus 2 as in Proposition 5.2. Here, we give a proof of Theorem 1.2, by using the previous results.

$\;$

{\it{Proof of Theorem 1.2}}. Let $X$ be a K3 surface whose N$\rm{\acute{e}}$ron-Severi lattice is a 2-elementary lattice with $\rho(X)=a=9$. Then there exists a non-symplectic involution $\theta$ which acts trivially on $S_X$, and, by Theorem 5.1, the set of the fixed points $X^{\theta}$ of it is a smooth curve of genus 2 which is ample as a divisor on $X$. Let $Y=X/\langle\theta\rangle$, and let $\pi:X\rightarrow Y$ be the natural quotient map. Since $\pi$ is the double cover branched along the smooth curve $\pi(X^{\theta})\in|-2K_Y|$, $Y$ is a DelPezzo surface given by the blow up at 8 points on general position in $\mathbb{P}^2$. Here, let $\varphi$ be the blow up at 8 points $P_i\;(1\leq i\leq 8)$ in $\mathbb{P}^2$, and let $\psi=\varphi\pi$. Moreover, let ${E_i}=\psi^{-1}(P_i)\;(1\leq i\leq 8)$ and let $B=\psi^{\ast}(l)$, for a line $l$ on $\mathbb{P}^2$. Then, since $S_X=\pi^{\ast}S_Y$, $S_X$ is generated by the classes of these curves. Since $X^{\theta}\in|-\pi^{\ast}K_Y|$, we have $X^{\theta}\sim 3B-\sum_{1\leq i\leq 8}E_i$. Hence, if we set
$$D_1=B-E_1-E_2,\;\;D_2=B-E_3-E_4,$$ 
$$D_3=B-E_5-E_6,\;\;D_4=B-E_7-E_8,$$
\noindent for each $1\leq i\leq4$, we have ${D_i}^2=-2$ and $X^{\theta}.D_i=2$. Since the member of $|D_i|\;(1\leq i\leq 4)$ is a $(-2)$-curve, by Proposition 5.2, $\mathcal{O}_X(D_i)\;(1\leq i\leq4)$ is ACM and initialized. Here, we construct families of simple ACM bundles of rank $n\geq2$, by using extensions of these line bundles.

$\;$

\noindent We consider the case where $n=2$. Note that, since 
$$|D_1-D_2|=|D_2-D_1|=\emptyset,$$
\noindent we have
$$\Hom(\mathcal{O}_X(D_1),\mathcal{O}_X(D_2))=\Hom(\mathcal{O}_X(D_2),\mathcal{O}_X(D_1))=0.$$
\noindent Hence, the assumption of Proposition 6.1 is satisfied, for the line bundles $\mathcal{O}_X(D_i)\;(i=1,2)$. Since
$$\Ext^2(\mathcal{O}_X(D_1),\mathcal{O}_X(D_2))\cong H^2(\mathcal{O}_X(D_2-D_1))=0,$$
\noindent we have
$$\dim\Ext^1(\mathcal{O}_X(D_1),\mathcal{O}_X(D_2))=-\chi(\mathcal{O}_X(D_2-D_1))=2.$$
\noindent In Proposition 6.1, if we set $U=\Ext^1(\mathcal{O}_X(D_1),\mathcal{O}_X(D_2))\backslash\{0\},$
\noindent we can construct a family parameterized by $U/\sim_w\cong\mathbb{P}^1$ of simple ACM bundles of rank 2. Here, we take $m$ distinct non-trivial extensions of $\mathcal{O}_X(D_1)$ by $\mathcal{O}_X(D_2)$
$$0\longrightarrow\mathcal{O}_X(D_2)\longrightarrow\mathcal{E}_i\longrightarrow\mathcal{O}_X(D_1)\longrightarrow0\;(1\leq i\leq m)\;\;\;\;\;\;(2)$$
\noindent such that any two extensions of them are not weak equivalent. Let 
$$\mathcal{E}:=\bigoplus_{1\leq i\leq m}\mathcal{E}_i.$$ 
\noindent Then we note that, by Proposition 6.1 (ii), we have $\Hom(\mathcal{E}_i,\mathcal{E}_j)=0$.

$\;$

\noindent We consider the case where $n=2m+1\;(m\geq1)$. Since 
$$|D_3-D_i|=|D_i-D_3|=\emptyset\;(i=1,2),$$ 
\noindent we have
$$\Hom(\mathcal{O}_X(D_3),\mathcal{O}_X(D_i))=\Hom(\mathcal{O}_X(D_i),\mathcal{O}_X(D_3))=0\;(i=1,2).$$
\noindent Hence, applying $\Hom(\mathcal{O}_X(D_3),-)$ and $\Hom(-,\mathcal{O}_X(D_3))$ to (2), we have
$$\Hom(\mathcal{E}_i,\mathcal{O}_X(D_3))=\Hom(\mathcal{O}_X(D_3),\mathcal{E}_i)=0\;(1\leq i\leq m).$$
\noindent Hence, the assumption of Proposition 6.1 is satisfied, for $m+1$ simple ACM bundles $\mathcal{E}_1,\cdots,\mathcal{E}_m$ and $\mathcal{O}_X(D_3)$. We consider the extension of $\mathcal{O}_X(D_3)$ by $\mathcal{E}$
$$0\longrightarrow\mathcal{E}\longrightarrow\mathcal{G}\longrightarrow\mathcal{O}_X(D_3)\longrightarrow0.$$
\noindent First of all, if we apply $\Hom(\mathcal{O}_X(D_3),-)$ to the exact sequence
$$0\longrightarrow \mathcal{O}(D_2)\longrightarrow\mathcal{E}_i\longrightarrow\mathcal{O}_X(D_1)\longrightarrow0,$$
\noindent we have
$$0\longrightarrow \Ext^1(\mathcal{O}_X(D_3),\mathcal{O}_X(D_2))\longrightarrow\Ext^1(\mathcal{O}_X(D_3),\mathcal{E}_i)\longrightarrow\Ext^1(\mathcal{O}_X(D_3),\mathcal{O}_X(D_1))\longrightarrow0,$$
\noindent and hence, we have $\dim\Ext^1(\mathcal{O}_X(D_3),\mathcal{E}_i)=4$. Therefore, in Proposition 6.1, if we set
$$U=\Ext^1(\mathcal{O}_X(D_3),\mathcal{E}_1)\backslash\{0\}\times\cdots\times\Ext^1(\mathcal{O}_X(D_3),\mathcal{E}_m)\backslash\{0\},$$
\noindent we get a family parameterized by $U/\sim_w\cong(\mathbb{P}^3)^m$ of simple ACM bundles of rank $2m+1$.

\pagebreak

\noindent We consider the case where $n=2m+2\;(m\geq1)$. Since 
$$|D_4-D_i|=|D_i-D_4|=\emptyset\;(1\leq i\leq3),$$
\noindent we have
$$\Hom(\mathcal{O}_X(D_3),\mathcal{O}_X(D_4))=\Hom(\mathcal{O}_X(D_4),\mathcal{O}_X(D_3))=0.$$
\noindent Applying $\Hom(\mathcal{O}_X(D_4),-)$ and $\Hom(-,\mathcal{O}_X(D_4))$ to (2), we have
$$\Hom(\mathcal{E}_i,\mathcal{O}_X(D_4))=\Hom(\mathcal{O}_X(D_4),\mathcal{E}_i)=0.$$
\noindent Here, we take an extension of $\mathcal{O}_X(D_3)$ by $\mathcal{E}$ which is a simple vector bundle
$$0\longrightarrow\mathcal{E}\longrightarrow\mathcal{G}\longrightarrow\mathcal{O}_X(D_3)\longrightarrow0.\;\;\;\;\;\;\;(3)$$
\noindent Since, if we apply $\Hom(\mathcal{O}_X(D_4),-)$ and $\Hom(-,\mathcal{O}_X(D_4))$ to (3), then we have
$$\Hom(\mathcal{O}_X(D_4),\mathcal{G})=\Hom(\mathcal{G},\mathcal{O}_X(D_4))=0,$$
\noindent for $\mathcal{O}_X(D_4)$ and $\mathcal{G}$, the assumption of Proposition 6.1 is satisfied. By the exact sequence
$$0\longrightarrow \Ext^1(\mathcal{O}_X(D_4),\mathcal{E})\longrightarrow\Ext^1(\mathcal{O}_X(D_4),\mathcal{G})\longrightarrow\Ext^1(\mathcal{O}_X(D_4),\mathcal{O}_X(D_3))\longrightarrow0,$$
\noindent we have
$$\dim\Ext^1(\mathcal{O}_X(D_4),\mathcal{G})=\dim\Ext^1(\mathcal{O}_X(D_4),\mathcal{O}_X(D_3))+\dim\Ext^1(\mathcal{O}_X(D_4),\mathcal{E}).$$
\noindent Applying $\Hom(\mathcal{O}_X(D_4),-)$ to (2), we have the exact sequence
$$0\longrightarrow \Ext^1(\mathcal{O}_X(D_4),\mathcal{O}_X(D_2))\longrightarrow\Ext^1(\mathcal{O}_X(D_4),\mathcal{E}_i)\longrightarrow\Ext^1(\mathcal{O}_X(D_4),\mathcal{O}_X(D_1))\longrightarrow0,$$
\noindent and hence, we have $\dim\Ext^1(\mathcal{O}_X(D_4),\mathcal{G})=4m+2\;(m\geq1)$. Therefore, in Proposition 6.1, if we set
$$U=\Ext^1(\mathcal{O}_X(D_4),\mathcal{G})\backslash\{0\},$$
\noindent we get a family parameterized by $U/\sim_w\cong\mathbb{P}^{4m+1}$ of simple ACM bundles of rank $2m+2$. 

The simple ACM bundles constructed by the above method are semistable. In fact, since the line bundles $\mathcal{O}_X(D_i)\;(1\leq i\leq 4)$ satisfy $H.D_i=H.D_j\;(1\leq i,j\leq 4)$ for $H\in|L|$, and ${D_i}^2=-2\;(1\leq i\leq 4)$, they have the same Hilbert polynomials. Since the simple ACM bundles constructed as above are given by the extensions of such  line bundles, by using induction and Lemma 6.1, we can easily see that they satisfy the semistability as in Definition 6.1 (however, they are not stable).$\hfill\square$

\smallskip

\smallskip

\noindent {\bf{Acknowledgements}}. The author would like to thank Prof. Konno. The author is partially supported by Grant-in-Aid for Scientific Research (25400039), Japan Society for the Promotion Science.


\begin{thebibliography}{Lambox}

\bibitem[B-P-W]{} W. Barth, C. Peters and A. van de Ven,
 \newblock{\em Compact complex surfaces}, 
  \newblock Springer, Berlin, (1984)

\bibitem[C-H]{} M. Casanellas - R. Hartshorne,
 \newblock{\em ACM bundles on cubic surfaces}, 
  \newblock J. Euro. Math. 13, 709-731 (2008)
  
\bibitem[C-K-M]{} E. Coskun, R. Kulkarni and Y. Mustopa, 
 \newblock{\em Pfaffian quartic surfaces and representations of Clifford algebras}, 
  \newblock Doc. Math. Vol. 17, 1003-1028 (2012)
  
\bibitem[C-P]{} C. Ciliberto and G. Pareschi, 
 \newblock{\em Pencils of minimal degree on curves on a K3 surface}, 
  \newblock J. reine angew. Math. 460 15-36 (1995)
  
%\bibitem[E-S-W]{} D. Eisenbud, F.-O. Schreyer, and J. Weyman,
 %\newblock{\em Resultants and Chow forms via exterior syzygies},
  %\newblock J. Amer. Math. Soc. 16, no. 3, 537-579 (2003)
  
\bibitem[Fa]{} D. Faenzi,
 \newblock{\em Rank 2 arithmetically Cohen-Macaulay bundles on a nonsingular cubic surface}, 
  \newblock J. Algebra 319 (1), 143-186 (2008)

\bibitem[Kn]{} H. Kn${\rm{\ddot{o}}}$rrer,
 \newblock{\em Cohen-Macaulay modules on hypersurface singularities I}, 
  \newblock Inv. Math. 88 (1), 153-164 (1987)
  
\bibitem[Mar]{} M.Maruyama,
 \newblock{\em Moduli of stable sheaves II}, 
  \newblock J. Math. Kyoto Univ. 3, 557-614 (1978)
  
\bibitem[M-M]{} S. Mori and S. Mukai,
 \newblock{\em The uniruledness of the moduli space of curves of genus 11}, 
  \newblock In Algebraic geometry, Lecture Notes in Math. 1016, Springer-Verlag, Berlin, 334-353 (1983). 
  
\bibitem[Ni]{} V.V. Nikulin, 
 \newblock{\em Factor groups of groups of automorphisms of hyperbolic forms with respect to subgroups generated by 2-reflections}, 
  \newblock J. Soviet Math. 22 1401-1476 (1983)
  
\bibitem[P-T]{} J. Pons-Llopis and F. Tonini,
 \newblock{\em ACM bundles on DelPezzo surfaces}, 
  \newblock Le Mathematiche 64, 177-211 (2009)
  
\bibitem[SD]{} B. Saint-Donat,
 \newblock{\em Projective Models of K3 surfaces}, 
  \newblock Amer. J. Math. Ann. 96, No.4, 602-639 (1974)
  
\bibitem[W1]{} K. Watanabe,
 \newblock{\em Donagi-Morrison's examples on 2-elementary K3 surfaces}, 
  \newblock Arch. Math. 98, 129-132 (2012)

\bibitem[W2]{} K. Watanabe,
 \newblock{\em The classification of ACM line bundles on quartic hypersurfaces in $\mathbb{P}^3$}, 
  \newblock Geometriae Dedicata, 1-8 (2014).
  
\end{thebibliography}
\end{document}